\documentclass{amsart}
\usepackage{amssymb,latexsym}
\theoremstyle{plain}
\newtheorem{theorem}{Theorem}

\newtheorem{lemma}{Lemma}

\newtheorem{proposition}{Proposition}
\theoremstyle{definition}

\theoremstyle{remark}

\numberwithin{equation}{section}
\newcommand{\R}{\mathbb R}

\newcommand{\C}{\mathbb C}

\begin{document}
%
%
% I slightly modify the short title
%
%
%
%
\title[Null--Controllability of 1--D parabolic equations]{Null--Controllability of one--dimensional parabolic equations}
%%%%%%%%%%%%%%%%%%
%Author information
\author{G. Alessandrini}
\address[G. Alessandrini]{Dipartimento di Matematica e Informatica\\Universit\`{a} degli Studi di Trieste\\Via Valerio, 12/b\\34127 Trieste, Italy.}
\email{alessang@univ.trieste.it}
\thanks{The first author is supported in part by MIUR, PRIN n. 2004011204.}
%%%%%%%%%%%%%%%%%%%%
\author{L. Escauriaza}
\address[L. Escauriaza]{Universidad del Pa{\'\i}s Vasco / Euskal Herriko
Unibertsitatea\\Dpto. de Matem\'aticas\\Apto. 644, 48080 Bilbao, Spain.}
\email{luis.escauriaza@ehu.es}
\thanks{The second author is supported  by MEC grant, MTM2004-03029.}
%\thanks{}
%%%%%%%%%%%%%%%%%%%%%%
%\keywords{Parabolic}
%\subjclass{Primary: subject; Secondary: subject}
%\date{}
%\dedicatory{}
%%%%%%%%%%%%%%
\begin{abstract}
We prove the interior null--controllability of one--dimensional
parabolic equations with time independent measurable coefficients.
\end{abstract}
\maketitle
%%%%%%%%%%%%%%%%%%%%%%

\section{Introduction}\label{S:1}
Let us consider the following one--dimensional heat equation with
variable coefficients
\begin{equation}\label{E: parab—lica}
\begin{cases}
\partial_x\left(a(x)\partial_xz\right)+b(x)\partial_xz+c(x)z-\rho(x)\partial_tz=f\chi_\omega, &\ 0<x<1\ ,\ 0<t<T\ ,\\
z(0,t)=z(1,t)=0, &\ 0\le t\le T\ ,\\
z(x,0)=z_0, & \ 0\le x\le 1\ .
\end{cases}
\end{equation}
Here, $z(x, t)$ is the state and $f(x,t)$ is an interior control that acts on the system over the open set $\omega\subset (0,1)$. The
coefficients $a$, $b$, $c$ and  $\rho$ are assumed to be measurable, bounded and for some $K \geq 1$,
\begin{equation}\label{E: condici—n1}
K^{-1}\le \rho (x)\le K,\quad  K^{-1}\le a(x)\le K,\ |b(x)|+|c(x)| \leq K \ , \ \text{a.e. in} \ [0,1]\ .
\end{equation}

For any given $z_0$ in $L^2(0, 1)$ and $f$ in $L^2(\omega\times [0, T])$, there is  only one solution $z$ to \eqref{E: parab—lica} in $C([0,
T];L^2(0, 1))\cap L^2(0,T; H^1_0(0,1))$. The goal of this paper is to analyze the interior null--controllability of \eqref{E: parab—lica}.
Specifically, we want to solve the following problem:

{\it Given $T > 0$ and $z_0$ in $L^2(0, 1)$, to find $f$ in $L^2(\omega\times [0, T])$ such that the corresponding solution $z$ to
 \eqref{E: parab—lica} satisfies, $z(\,\cdot\,, T) \equiv 0$ in $(0, 1)$.}

In \cite{iy}, it is shown that the system \eqref{E: parab—lica} is
null--controllable at any positive time, when the coefficients $a$
and $\rho$ are Lipschitz in $[0,1]$. In this reference, the proof
of null--controllability is based on an appropriate observability
inequality for the adjoint system and it is implied by a global
Carleman estimate. When the coefficients  are smooth, the
observability inequality can be proved introducing Fourier series
and using high frequency asymptotic formulae for the eigenvalues
of the corresponding St\"urm-Liouville problem and classical
results on the sums of real exponentials, see \cite{loz}. In
\cite{fz2} adopting the approach introduced by D.L. Russel in
\cite{r} (the null controllability of the wave operator at large
times implies the null controllability of the heat equation at all
times) it is shown that the system \eqref{E: parab—lica} is null
controllable, when $a$ and  $\rho$ have bounded variation in
$[0,1]$.

The main result in this paper is the following.
\begin{theorem}\label{T: 1}
Assume that the coefficients $a$, $b$, $c$ and $\rho$ are bounded measurable and satisfy \eqref{E: condici—n1}. Then, \eqref{E: parab—lica} is
null--controllable at time $T$, for all $T>0$ and with controls $f$ in $L^2(0,T; H^1_0(\omega))$.
\end{theorem}

To prove this result we proceed in the the following way. First, a change of variables shows that the internal controllability of the system
\eqref{E: parab—lica} is equivalent to the same question for a system
%
%
%I think we can maintain the interval (0,1)
%
%
%
\begin{equation}\label{E: parab—lica2}
\begin{cases}
\partial_{x}^2z-\rho(x)\partial_tz=f\chi_\omega, &\ 0<x<1\ ,\ 0<t<T\ ,\\
z(0,t)=z(1,t)=0, &\ 0\le t\le T\ ,\\
z(x,0)=z_0, & \ 0\le x\le 1\ ,
\end{cases}
\end{equation}
where $\rho$ is a new measurable function satisfying \eqref{E: condici—n1}, for some new constant only depending on $K$, which we shall continue
to denote by $K$.

Then, if we denote by $\delta$ the inradius of the open set $\omega$, that is
\begin{equation} \label{E: inradius}
\delta =\sup \{ r> 0 |\  \exists \xi \in \omega  \ , \ (\xi -r,
\xi+r) \subset \omega \} \ ,
\end{equation}
%
%
% I have replaced  \omega _k with \lambda_k , in order not to mix the eigenvalues with the control set \omega
%
% also note that such eigenvalues are all simple
%
%
and if $e_1, e_2,\dots , e_n \dots$ and $0<\lambda_1^2<\lambda_2^2<\dots<\lambda_m^2\dots$ are respectively the eigenfunctions and eigenvalues
of the St\"urm-Liouville problem
\begin{equation}\label{E:sturm}
\begin{cases}
e''+\rho(x)\lambda^2e=0\ , \ 0<x<1 \ ,\\
e(0)=e(1)=0\ ,\\
\end{cases}
\end{equation}
we prove the following:
\begin{theorem}\label{T: 2}
Assume that the coefficient $\rho$ satisfies \eqref{E: condici—n1}. Then, there is a constant $N$, which depends on $K$ and on $\delta$
 such that the inequality
\begin{equation}\label{E: propiedad}
\sum_{\lambda_k\le\mu}a_k^2\le Ne^{N\mu}\int_{\lambda}|\sum_{\lambda_k\le\mu}a_ke_k |^2\,dx\ ,
\end{equation}
holds for all $\mu\ge 1$ and all sequences $\{a_k\}$.
\end{theorem}
The analog of this inequality for the eigenfunctions of the Laplace-Beltrami operator  on a compact and {\it smooth} Riemannian manifold with a
possibly nonempty boundary was proved in \cite{lr}. There, G. Lebeau and L. Robbiano showed that it implies the interior null-controllability of
the heat equation over the manifold by giving an explicit construction of the control function $f$ (See \cite[\S 5]{lz} for a more simplified
presentation).

The arguments in \cite{lr} show that the same iterative method of construction of the control function $f$ given in \cite{lr} works for the
system \eqref{E: parab—lica2}, when Theorem \ref{T: 2} holds. Thus, Theorem \ref{T:  1} follows from Theorem \ref{T: 2}.

To prove Theorem \ref{T: 2} we start by following the arguments in \cite{lr}. In particular, given $\mu\ge 1$ and a sequence of real numbers
$a_1,a_2,\dots ,a_n\dots$, we set
\begin{equation*}
u(x,y)=\sum_{\lambda_k\le\mu}a_ke_k(x)\cosh{\left(\lambda_ky\right)}\ .
\end{equation*}
This function satisfies
\begin{equation}\label{E: eliptico2}
\begin{cases}
\partial^2_{x}u+\partial_y\left (\rho(x)\partial_{y}u\right)=0, &\ 0<x<1\ ,\ y\in\R\ ,\\
u(0,y)=u(1,y)=0, &\ y\in\R\ ,\\
\partial_yu(x,0)=0, &\ 0<x<1\ ,\\
\end{cases}
\end{equation}
and the proof of Theorem \ref{T: 2} is a consequence of a {\it quantification} of the following  {\it qualitative} result of unique continuation
from the boundary:

\vspace{0.1cm}
{\it Assume that $u$ satisfies \eqref{E: eliptico2} and $u(x,0)\equiv 0$, when $x$ is in $\omega\subset (0,1)$.Then, $u\equiv 0$ in $[0,1]\times\R$.}
\vspace{0.1cm}

In \cite{lr}, the one dimensional interval $[0,1]$ is replaced by a compact and {\it smooth} manifold $M$, $\partial^2_{x}$  by the
corresponding Laplace-Beltrami operator on $M$ and the authors work out  the {\it quantification} of a similar {\it qualitative} property of
boundary unique continuation for  the elliptic operator,  $\triangle+\partial^2_{y}$, where $\triangle$ is the Laplace-Beltrami operator on $M$.
To carry out this quantification they use two Carleman inequalities. Those methods require that the elliptic operator involved has Lipschitz
second order coefficients and so, they can not be applied to the elliptic operator in \eqref{E: eliptico2}, which has  measurable coefficients.

On the other hand, if $\sigma$ is a $2\times 2$ symmetric and measurable matrix in the plane verifying  the ellipticity condition
\begin{equation}\label{E:condicion2}
K^{-1}|\xi|^2\le \sigma(x,y)\xi\cdot\xi\le K|\xi|^2,\ \text{when}\ (x,y)\ \text{and}\ \xi\in\R^2\ ,
\end{equation}
the weak solutions of the equation
\begin{equation}\label{E:eliptico2}
\nabla\cdot\left(\sigma (x,y)\nabla u\right)=0\ ,
\end{equation}
 satisfy the strong unique continuation property:

\vspace{0.1cm} {\it If a $W^{1,2}_{loc}$-solution of  \eqref{E:eliptico2} on a connected open set $\Omega$ has a zero of infinite order at an
interior point, then it must be zero.} \vspace{0.1cm}

See \cite{am}. This qualitative result of strong unique continuation for uniformly elliptic  equations in two independent variables is based on
the connection between the solutions of these equations and the theory of quasiregular mappings \cite{bn} and on  the so-called Ahlfors-Bers
representation \cite{ab} of such mappings. Here, we describe some quantifications  of this qualitative result and apply them to prove the
null--controllability property. In particular, a \ \lq\lq Hadamard's three circle theorem\rq\rq\ , Proposition \ref{P:threesph}, and a \ \lq\lq
doubling\rq\rq \ type property, Proposition \ref{P:doubling2}, adapted to the solutions of \eqref{E:eliptico2}.

In section \ref{S:2} we  recall the results we need from the theory of quasiregular mappings and  prove the adapted Hadamard's three circle
theorem and doubling property. In section \ref{S:3} we show how to apply them to prove Theorem \ref{T: 2}, also using an estimate of
continuation from Cauchy data for solutions of \eqref{E:eliptico2} Lemma \ref{L:1}, which we adapt from \cite{ar}. It may be worth noting that
the approach used for the proof of Lemma \ref{L:1}, is based on a variation on the classical principle of majorization by harmonic measure,
\cite[Chapter VIII, \S 1, p. 301]{tsu}, which in turn has its roots in arguments due to Carleman \cite[p. 3--4]{carl}.

\section{Quantitative estimates of unique continuation with discontinuous coefficients}\label{S:2}
%
%
% I have shortened the part of introduction to qc mappings
%
%
Throughout the paper, $z=x+iy$, $\Omega$ is a simply connected open set in the plane, $B_r$ a circle of radius $r$ centered at the origin, and
\begin{equation*}
\partial_{\overline z}f=\tfrac 12\left(\partial_xf+i\partial_yf\right)\ ,\ \partial_{z}f=\tfrac 12\left(\partial_xf-i\partial_yf\right)\ .
\end{equation*}
We shall denote by $C$ constants only depending on $K$, whereas by $N$ we shall denote constants only depending on $K$ and $\delta$.

When $u \in W^{1,2}_{loc}(\Omega)$ is a weak solution to  \eqref {E:eliptico2}, and $\sigma$ satisfies \eqref{E:condicion2}
 we can associate in a natural fashion, which generalizes the harmonic conjugate, a new function, the so called {\it stream} function $v$, which satisfies
\begin{equation}\label{E:sistema}
\nabla v = J \sigma \nabla u
\end{equation}
almost everywhere in $\Omega$ and is a weak solution to
\begin{equation}\label{E:eliptico3}
\nabla\cdot\left(\tfrac{\sigma}{\det{\sigma}}\nabla v\right)=0,\  \text{in}\ \Omega\ .
\end{equation}
Here $J$ denotes the matrix representing a $90^{\circ}$ rotation in the plane
\begin{equation*}
J = \left( \begin{array}{cc}
   0&-1\\
   1& 0\\
 \end{array}\right) \ .
\end{equation*}
 Moreover, letting $f = u + iv$, we have $f\in W^{1,2}_{loc}$ and
 satisfies
\begin{equation}\label{E:quasi-regular}
\partial_{\overline z}f=\mu\partial_zf+\nu\overline{\partial_{z}f},\ \text{almost everywhere in $\Omega$}\ ,
\end{equation}
where the complex valued functions $\mu$ and $\nu$ can be explicitly expressed in terms of $\sigma$ , see \cite{am}, and verify
\begin{equation}\label{E:quesi-regular condition}
|\mu|+|\nu|\le \tfrac{K-1}{K+1}<1,\ \text{almost everywhere in $\Omega$}
\end{equation}
That is, $f$ is a $K$-quasiregular mapping.

To give an idea of why these results hold,  observe that the
vector field
\begin{equation*}
J \sigma \nabla u
\end{equation*}
is, in the weak sense, curl-free in $\Omega$. To verify that $v$
is a $W^{1,2}_{loc}$-solution of \eqref{E:eliptico3}, observe
that, from \eqref{E:sistema}, one obtains that the vector field
\begin{equation*}
\tfrac{\sigma}{\det{\sigma}}\nabla v = J {\sigma}^{-1}J^t\nabla v
\end{equation*}
is, in the weak sense, divergence-free in $\Omega$.

By the Ahlfors-Bers representation \cite{ab} (see also \cite{bn}
and  \cite[Chapter II.6, pp. 258--259]{bjs}), any K-quasiregular
mapping $f$ in $B_1$ can be written as
\begin{equation*}
f=F\circ\chi\ ,
\end{equation*}
where $F$ is holomorphic in $B_1$ and $\zeta=\chi(z)$ is a $K$-quasiconformal homeomorphism from $B_1$ onto $B_1$, which verifies, $\chi(0)=0$,
$\chi(1)=1$,
\begin{equation}\label{E:Hšlder}
C^{-1}|z_1-z_2|^{\frac 1\alpha}\le |\chi(z_1)-\chi(z_2)|\le
C|z_1-z_2|^\alpha\ ,\text{\ when $z_1, z_2\in B_1$}
\end{equation}
for some $0<\alpha<1$ and $C\ge 1$ depending only on $K$.

We now recall the Hadamard's three-circle theorem \cite{m}.

\begin{theorem}\label{T: 3}
 Let $F$ be a holomorphic function of a complex variable in the ball $B_{r_2}$ and $M(r)=\max_{B_r}{|F|}$. Then,
 the following is valid for $0<r_1\le r \le r_2$,
\begin{equation*}\label{E:Hadamard}
\log {M(r)}\le \frac{\log{\frac
{r_2}{r}}}{\log{\frac{r_2}{r_1}}}\log {M(r_1)}+\frac{\log{\frac
{r}{r_1}}} {\log{\frac{r_2}{r_1}}}\log {M(r_2)} \ .
\end{equation*}
\end{theorem}
The meaning of this inequality is that  $\log{M(r)}$ is a convex function of the variable $\log{r}$.

Let $u \in W^{1,2}_{loc}(B_R)$ be a weak solution to
\eqref{E:eliptico2} and let $f:B_R\longrightarrow\C$ be the
associated $K$-quasiregular mapping. Rescaling \eqref{E:Hšlder} we
have that $f=F\circ\chi$, where $F$ is holomorphic in $B_R$ and
$\chi : B_R\longrightarrow B_R$ is a $K$-quasiconformal
homeomorphism, which verifies
\begin{equation}\label{E:dilataci—n2}
RC^{-1}|\tfrac zR|^{\frac 1\alpha}\le |\chi(z)|\le RC|\tfrac
zR|^\alpha\
\end{equation}
where $C$ is the same as in \eqref{E:Hšlder}.
%
%
% I have changed \xi with \zeta , which is easily understood as a complex quantity, and accordingly changed \varsigma with \xi
%
%
%

 Define
\begin{equation}\label{E:twistedballs}
 \mathcal B_r=\{z\in B_R : |\chi(z)|<r\}
\end{equation}
and
\begin{equation}\label{E:twistedmax}
m(r)=\max_{\mathcal B_r}|f(z)|\ , \text{ when } r<R \ .
\end{equation}
Then, through the change of coordinates, $\zeta=\chi(z)$, the
Hadamard's three circle theorem takes the form:  the function
$\log{m(r)}$ is a convex function of $\log r$,
\begin{equation}\label{E:Hadamard2}
\log {m(r)}\le \frac{\log{\frac
{r_2}{r}}}{\log{\frac{r_2}{r_1}}}\log {m(r_1)}+\frac{\log{\frac
{r}{r_1}}}{\log{\frac{r_2}{r_1}}}\log {m(r_2)}\ ,\ \text{when
$0<r_1\le r\le r_2 < R$,}
\end{equation}
and the sets $\mathcal B_r$, almost look like balls. In particular,
\begin{equation}\label{E:lasbolas}
\mathcal B_R=B_R\ \ ,\ B_{R\left(\frac r{CR}\right)^{\frac
1\alpha}}\subset\mathcal B_r\subset B_{R\left(\frac
{Cr}R\right)^{\alpha}}\ ,\ \text{\ when $r<R$}
\end{equation}
and $C$ is the same constant appearing in \eqref{E:dilataci—n2}.
Note incidentally that \eqref{E:Hadamard2} implies a weak unique
continuation property, that is, if $m(r_1)=0$ for some small
$r_1$, then $m(r)=0$ for all $r<R$.

On the other hand, the difference quotients of convex functions are nondecreasing functions of their arguments. This implies that, if $f$ is not
identically zero,
\begin{equation*}
\frac{\log{m(\frac r2)}-\log{m(\frac r4)}}{\log {\frac r2}-\log{\frac r4}}\le \frac{\log{m(\frac R2)}-\log{m(\frac R4)}} {\log {\frac
R2}-\log{\frac R4}}\ ,\ \text{\ when}\ r\le R\ ,
\end{equation*}
and thus,
\begin{equation}\label{E:doubling}
\frac{m(\frac r2)}{m(\frac r4)}\le \frac {m(\frac R2)}{m(\frac R4)}\ ,\ \text{\ when}\ r<R\ .
\end{equation}
We may prescribe that the the stream function $v$ of $u$ satisfies $v(0)=0$. We have that $F=u+iv$ is holomorphic in the $\zeta=\xi+i\eta$
coordinates in $B_R$, hence, solving the Cauchy-Riemann equations,
\begin{equation*}
v(\xi,\eta)=\int_0^\eta u_\xi(\xi,s)\,ds-\int_0^\xi u_\eta(t,0)\,dt\ , \text{\ in $B_R$}\ .
\end{equation*}
This formula and interior estimates for harmonic functions
\cite{gt} show that in the $\zeta$-coordinates we have,
\begin{equation*}
\|u\|_{L^\infty(B_r)}\le \max_{B_r}|F(\zeta)|\le
C\|u\|_{L^\infty(B_{2r})}\ ,\ \text{when}\ r\le \frac R2\
\end{equation*}
where $C>0$ is an absolute constant. In the $z$-coordinates, the
last inequality reads as
\begin{equation}\label{E:comparaci—n}
\|u\|_{L^\infty(\mathcal B_r)}\le \max_{\mathcal B_r}|f(z)|\le C\|
u\|_{L^\infty(\mathcal B_{2r})}\ ,\ \text{when}\ r\le \frac R2\ ,
\end{equation}
and from \eqref{E:comparaci—n}, \eqref{E:Hadamard2} we obtain:

\begin{proposition} \label{P:threesph}
Let $u \in W^{1,2}_{loc}(B_R)$ be a weak solution to
\eqref{E:eliptico2} and let $\mathcal B_r \ , \ 0<r  \le R$ be the
open sets introduced in \eqref{E:twistedballs}, then we have
\begin{equation}\label{E:treesph}
\|u\|_{L^\infty(\mathcal B_{\frac r2})}\le
C\|u\|_{L^\infty(\mathcal B_{r_1})}^\theta\|u\|_{L^\infty(\mathcal
B_{r_2})}^{1-\theta}\ , \text{\ when $r_1\le r\le r_2 < R$ ,
$\theta=\frac{\log{\frac{r_2}{r}}}{\log{\frac{r_2}{r_1}}}$} \ .
\end{equation}
\end{proposition}

And also, from  \eqref{E:doubling}:

\begin{proposition} \label{P:doubling2}
Let $u$  and $\mathcal B_r $ be as above. If $u$ is not
identically zero, then we have
\begin{equation}\label{E:doubling2}
\frac{\|u\|_{L^\infty(\mathcal B_{r})}}{\|u\|_{L^\infty(\mathcal
B_{\frac r2})}}\le C\frac{\|u\|_{L^\infty(\mathcal
B_{R})}}{\|u\|_{L^\infty(\mathcal B_{\frac R4})}}\ , \text{\ when
$r\le R$}\ .
\end{equation}
\end{proposition}

These are respectively a Hadamard's three circle theorem
%or three sphere's inequality
and a doubling property
adapted to the solution $u$ through the family of \ \lq\lq balls\rq\rq\ $\mathcal B_r$.

\section{Proof of Theorem \ref{T: 1}}\label{S:3}
First let us note that, possibly replacing $z$ in \eqref{E: parab—lica} with $e^{-2Kt}z$, we may assume that $c$ is nonpositive. Introducing
$B(x)=\int_0^x\frac{b(s)}{a(s)}\,ds$ , we observe that \eqref{E: parab—lica} can be rewritten as
\begin{equation*}
e^{-B(x)}\partial_x\left(a(x)e^{B(x)}\partial_xz\right)+c(x)z-\rho(x)\partial_tz=f\chi_\omega\ ,\ \text{\ in $(0,1)\times (0,T]$}\ .
\end{equation*}
The solution $w$ to
\begin{equation*}
\begin{cases}
e^{-B(x)}\frac {d}{dx}\left(a(x)e^{B(x)}\frac{dw}{dx}\right)+c(x)w=0\ ,\\
w(0)=w(1)=1\ ,
\end{cases}
\end{equation*}
verifies, $0<w(x)\le 1$ in $[0,1]$ and replacing $z$ with the new dependent variable $\tilde z =z/w$, which we denote again $z$, we have
\begin{equation*}
e^{-B(x)}\partial_x\left(a(x)w^2(x)e^{B(x)}\partial_xz\right)-\rho(x)w^2(x)\partial_tz=w(x)f\chi_\omega\ ,\ \text{\ in $(0,1)\times (0,T]$}\ .
\end{equation*}
Setting
\begin{equation*}
L=\int_0^1a^{-1}(s)w^{-2}(s)e^{-B(s)}\,ds \ , \
y=\tfrac{1}{L}\int_0^xa^{-1}(s)w^{-2}(s)e^{-B(s)}\,ds
\end{equation*}
and writing $z(x,t)=\tilde z(y,t)$,  the new function $\tilde z$,
which again we rename $z$, is a solution  of a system of the form
\eqref{E: parab—lica2}.

Considering the associated St\"urm-Liouville problem
\eqref{E:sturm}, we extend the eigenfunctions $e_j$, $j\ge 1$, to
$[-1,1]$ by an odd reflection in $0$, similarly we extend $\rho$
by an even reflection in $0$. Next, we continue these new
functions to all of $\R$ as periodic functions of period $2$. The
extended $\rho$ verifies \eqref{E: condici—n1}, $e_j\in
C^{1,1}(\R)$ and $e_j''+\rho(x)\lambda_j^2e_j=0$, almost
everywhere in $\R$.

Being the change of variable $y=y(x)$ bi--Lipschitz, with Lipschitz constants which only depend on $K$, the open set $\omega$ is transformed
into a new open subset of $(0,1)$ whose inradius is comparable to $\delta$. We continue to denote the transformed set and its inradius by
$\omega$ and $\delta$, respectively. Also we can assume, up to a translation along the real line, $(-\delta,\delta)\subset\omega\subset (-1,1)$.

Given $\mu\ge 1$ and a sequence of real numbers $a_1, a_2,\dots a_n,\dots$, the function
\begin{equation*}
u(x,y)=\sum_{\lambda_k\le\mu}a_ke_k(x)\cosh{\left(\lambda_ky\right)}\ .
\end{equation*}
verifies
\begin{equation}
\begin{cases}
\partial^2_{x}u+\partial_y\left (\rho(x)\partial_{y}u\right)=0\ ,\ \text{\ in $\R^2$}\ ,\\
\partial_yu(x,0)=0\ ,\ \text{\ in $\R$}
\end{cases}\end{equation}
and its stream function $v$ can be chosen so that, $v\in W^{1,2}_{loc}(\R^2)$,
\begin{equation*}
\begin{cases}
\partial_xv=-\rho(x)\partial_yu\ ,\\
\partial_yv=\partial_xu\ ,
\end{cases}\quad \text{and}\quad
\begin{cases}
\partial_x\left(\frac 1{\rho(x)}\partial_xv\right)+\partial^2_yv=0\ ,\ \text{\ in $\R^2$,}\\
v(x,0)=0,\  \text{\ in $\R$}\ .
\end{cases}
\end{equation*}
Let $f=u+iv$, consider the family of \ \lq\lq balls\rq\rq\
$\mathcal B_r$ associated to $f$ in section \ref{S:2} at scale $R$
and choose $R=2\left(4C\right)^{\frac 1\alpha}$, where $\alpha$
and $C$ are the constants in \eqref{E:dilataci—n2} and
\eqref{E:lasbolas}. With this choice, $\mathcal B_{\frac
R4}\supset B_2$. The interior bounds for subsolutions of elliptic
equations \cite{gt} give
\begin{equation*}
\|u\|_{L^\infty(B_R)}\le \tfrac CR\|u\|_{L^2(B_{2R})}\ .
\end{equation*}
These and the orthogonality of the eigenfunctions $e_j$, $j\ge 1$ imply that
\begin{equation}\label{E:doubling constant}
\frac{\|u\|_{L^\infty(\mathcal B_{R})}}{\|u\|_{L^\infty(\mathcal
B_{\frac R4})}}\le e^{C\mu}\ .
\end{equation}
An iteration of \eqref{E:doubling2} and \eqref{E:doubling constant} give
\begin{equation*}
\|u\|_{L^\infty(\mathcal B_{R})}\le
e^{Ck\mu}\|u\|_{L^\infty(\mathcal B_{\frac R{2^k}})}\ , \text{\
when $k\ge 1$}
\end{equation*}
and from \eqref{E:lasbolas}, there is $k=k(\delta ,K)$ such that, $\mathcal B_{\frac R{2^k}}\subset B_{\frac \delta 2}$. Thus,
\begin{equation}\label{E:doubling4}
\|u\|_{L^\infty(B_1)}\le e^{N\mu}\|u\|_{L^\infty(B_{\frac\delta 2})}\ .
\end{equation}

The following inequality, which is an estimate on the continuation
from Cauchy data, holds.
\begin{lemma}\label{L:1} There are constants $0<\theta<1$ and $C>0$, only depending on $K$, such that the inequality
\begin{equation*}
\|u\|_{L^\infty(B_{\frac r2})}\le Cr^{-\frac\theta 2}\|u(\,\cdot\, ,0)\|_{L^2(-r,r)}^\theta\|u\|_{L^\infty(B_{4r})}^{1-\theta}
\end{equation*}
holds, when $r\le 1$.
\end{lemma}
The Lemma and \eqref{E:doubling4} give
\begin{equation*}
\|u\|_{L^\infty(B_1)}\le N e^{N\mu}\|u(\,\cdot\, ,0)\|_{L^2(-\delta,\delta)}\ ,
\end{equation*}
which proves Theorem \ref{T: 2}.

\begin{proof}[Proof of Lemma \ref{L:1}]
%
%
% I have slightly modified this proof, so to avoid the use of the special fact that \u and \v are  C^1
% which would not hold for a general \sigma
%
%
The proof is essentially contained in  \cite[Theorem 4.5]{ar}. For the sake of completeness we summarize here the argument.

Recalling that $f=u+iv$ is analytic in the $\zeta$-variable, we
have that  $\log{|f(\zeta)|}$ is subharmonic in $B_R$ and
consequently one can verify that  $\log{|f(z)|}$ is a subsolution
for an elliptic operator in divergence form $E$ with a matrix of
coefficients verifying \eqref{E:condicion2}. For $r>0$, let $w$ be
the solution to
\begin{equation}
\begin{cases}
Ew=0\ ,\ \text{\ in $B_r^+$},\\
w=1\ ,\ \text{\ in $(-r,r)$},\\
w=0\ ,\ \text{\ in $\partial B_r^+\setminus (-r,r)$}\ .
\end{cases}
\end{equation}
On $\partial B_r^+$ we have
\begin{equation*}
\log{|f|}\le w\log{\|u(\,\cdot\, ,0)\|_{L^\infty (-r,r)}}+(1-w)\log{\|f\|_{L^\infty(B_{r})}}
\end{equation*}
and the maximum principle implies that the same inequality also holds in $B_r^+$. The H\"older continuity at the boundary of $w$ and the
Harnack's inequality \cite{gt} show that there is $\eta \in (0,1)$, which only depends on $K$,
 such that $w(z)\ge \eta$ in $B_{\frac r2}^+$.

Using $v(0)=0$, \eqref{E:sistema} and interior bounds for elliptic equations \cite{gt}, we have
\begin{equation*}
\|v\|_{L^\infty(B_{r})} \le C \|\nabla v\|_{L^2(B_{2r})} \le C \|\nabla u\|_{L^2(B_{2r})} \le C \|u\|_{L^\infty(B_{4r})} \ .
\end{equation*}

These imply
\begin{equation}\label{E:cauchy}
\|u\|_{L^\infty(B_{\frac r2})}\le C\|u(\,\cdot\, ,0)\|_{L^\infty(-r,r)}^\eta\|u\|_{L^\infty(B_{4r})}^{1-\eta}\ .
\end{equation}

For every $\alpha \in (0,1]$, we have the interpolation inequality
\begin{equation}\label{E:interp}
\|\varphi\|_{L^\infty(-r,r)}\le C\left(\|\varphi\|_{L^2(-r.r)}^{\beta}|\varphi |_{C^{\alpha}(-r,r)}^{1-\beta} + r^{-\frac
12}\|\varphi\|_{L^2(-r.r)} \right) \ ,
\end{equation}
where $\beta = \tfrac{2\alpha}{1+2\alpha}$, $C>0$ only depends on $\alpha$ and $|\varphi |_{C^{\alpha}(-r,r)}$ denotes the standard $C^{\alpha}$
seminorm. Next, we use the interior H\"older bound for $u$, \cite{gt},
\begin{equation}\label{E:interior}
|u|_{C^{\alpha}(B_{r})} \le C r^{-\alpha} \|u\|_{L^\infty(B_{4r})} \ ,
\end{equation}
with $C>0$ and $\alpha \in (0,1]$ only depending on $K$. Combining \eqref{E:cauchy} with \eqref{E:interp} and \eqref{E:interior}, we obtain the
thesis with $\theta= \beta \eta$.

The interpolation inequality \eqref{E:interp} can be proved essentially along the same lines as the interpolation inequalities in \cite[\S 6.8]{gt}.

\end{proof}
  %%%%%%%%%%%%%%%%%
\newpage
%%%%%%%%%%%%%%%%%%%%%%

%%%%%%%%%%%%%%%%%%%%%%%%%%

\begin{thebibliography}{99}

\bibitem{ab} L. Ahlfors, L. Bers, \emph{Riemann's mapping theorem for variable metrics,} Ann. Math. \textbf{72} (1960), 265--296.

%\bibitem{ahl} L. Ahlfors, \emph{Lectures on Quasiconformal Mappings,} Van Nostrand, Priceton N.J., 1966.

\bibitem{am} G. Alessandrini, R. Magnanini, \emph{Elliptic equations in divergence form, geometric critical
oints of solutions and Stekloff eigenfunctions,} SIAM J. Math. Anal., \textbf{25} n. 5 (1994), 1259--1268.

\bibitem{ar} G. Alessandrini, L. Rondi, \emph{Stable determination of a crack in a planar inhomogeneous conductor,} SIAM
J. Math. Anal. \textbf{30} n. 2 (1998), 326--340.

\bibitem{bjs} L. Bers, F. John, M. Schechter, \emph{Partial Differential Equations,} Interscience, New York, 1964.

\bibitem{bn} L. Bers, L. Nirenberg, \emph{On a representation theorem for linear elliptic systems with discontinuous
coefficients and applications ,} in Convegno Internazionale sulle Equazioni alle Derivate Parziali, Cremonese, Roma, 1955, 111--138.

\bibitem{az} C. Castro, E. Zuazua \emph{Concentration and lack of observability of waves in highly heterogeneous
media,} Archive Rational Mechanics and Analysis \textbf{164} n. 1 (2002), 39--72.

\bibitem{carl} T. Carleman, \emph{Les Fonctions Quasi Analytiques,} Gauthier--Villars, Paris 1926.

\bibitem{fz1} E. Fern\'andez-Cara, E. Zuazua, \emph{The cost of approximate controllability for heat equations:
The linear case,} Advances Diff. Eqs. \textbf{5} (4-6) (2000), 465--514.

\bibitem{fz2} \bysame, \emph{On the null controllability of the one-dimensional heat equation with BV coefficients,} Computational
and Applied Mathematics \textbf{21} (1) (2002), 167--190.

\bibitem{fi} A. V. Fursikov, O. Yu. Imanuvilov, \emph{Controllability of evolution equations,} Lecture Notes Series \textbf{34},
Research Institute of Mathematics, Global Analysis Research Center, Seoul National University (1996).

\bibitem{gt} D. Gilbarg, N.S. Trudinger, \emph{Elliptic Partial Differential Equations of Second Order,} 2nd ed., Springer-Verlag,
Berlin-Heildeberg-New York-Tokyo, 1983.

\bibitem{iy} O. Yu Imanuvilov, M. Yamamoto, \emph{Carleman estimate for a parabolic equation in Sobolev
spaces of negative order and its applications,} in Control of Nonlinear Distributed Parameter
Systems, G. Chen et al. eds., Marcel-Dekker  (2000), 113--137.

\bibitem{lo} E.M. Landis, O.A. Oleinik, \emph{Generalized analyticity and some related
properties of solutions of elliptic and parabolic equations,} Russian Math. Surv.
\textbf{29} (1974), 195--212.

\bibitem{lr} G. Lebeau, L. Robbiano, \emph{Contr\^ole exact de l'\'equation de la chaleur,}  Commun. Partial Differ. Equ.
\textbf{20} (1995), 335-356.

\bibitem{lz} G. Lebeau, E. Zuazua, \emph{Null controllability of a system of linear thermoelasticity,}  Archive for Rational Mechanics and Analysis
\textbf{141} (4) (1998), 297--329.

\bibitem{l} F.H. Lin, \emph{A uniqueness theorem for parabolic equations,} Comm. Pure
Appl. Math. \textbf{42} (1988), 125--136.

\bibitem{loz} A. L\'opez, E. Zuazua, \emph{Uniform null-controllability for the one-dimensional heat equation with rapidly oscillating periodic density,} Ann. I.H.P. - Analyse non lin\'eaire \textbf{19}, 5 (2002), 543--580.

\bibitem{m} A.I. Markushevich, \emph{Theory of Functions of a Complex Variable,} Prentice Hall, Englewood Cliffs, NJ, 1965.

\bibitem{r} D.L. Russel, \emph{A unified boundary controllability theory for hyperbolic and parabolic partial differential equations,} Studies in Appl. Math. \textbf{52} (1973), 189--221.

\bibitem{tsu} M. Tsuji, \emph{Potential Theory in Modern Function Theory,} Maruzen, Tokyo 1959.

\bibitem{v} I.N. Vekua, \emph{Generalized Analytic Functions,} Pergamon, Oxford 1962.

\end{thebibliography}
\end{document}